\documentstyle[12pt]{article}
\newcommand{\R}{\mbox{$I \kern -4pt R$}}             % simbolo de Reales
\oddsidemargin= - 0.5cm
\begin {document}
\baselineskip 30pt
\title{On pressure boundary conditions for thermoconvective problems} 
 \author{H. Herrero$\star$ and A. M. Mancho$_\dagger$ \footnote{Present address: School of Mathematics,  University of Bristol. University Walk, Bristol, BS8 1TW. United Kingdom.}\\
             {$\star$ \small Departamento de Matem\'aticas,} 
 	{\small Facultad de Ciencias Qu\'{\i}micas}\\
        {\small Universidad de Castilla-La Mancha,}
        {\small 13071 Ciudad Real, Spain.}\\
{\small e-mail: hherrero@qui-cr.uclm.es, phone: 
34 926 295412, fax: 34 926 295318}\\
{$\dagger$ \small Centro de Astrobiolog\'{\i}a, (Associate Member of NASA Astrobiology Institute)}\\
{\small CSIC-INTA, Ctra. Ajalvir, Km. 4, Torrej\'on 
de Ardoz, 28850 Madrid, Spain.}}
\maketitle


\noindent
\begin{abstract}
Solving numerically hydrodynamical problems of incompressible fluids raises 
the question of handling first order derivatives (those of pressure) in a closed
container and determining its boundary conditions. 
A way to avoid the first point is to derive a Poisson equation for pressure,
although the problem of  taking the 
right boundary conditions still remains. 
To remove this problem another formulation of the problem has been
used  consisting of projecting the master equations into 
the space of divergence free velocity fields, so pressure is eliminated
from the equations. This technique 
raises the order of the differential equations and additional 
boundary conditions may be required. 
High order derivatives are sometimes 
 troublesome, specially in cylindrical coordinates due to
 the singularity at the origin, so for these problems a low order formulation
 is very convenient.   
We research several pressure boundary conditions for the primitive
variables formulation of thermoconvective problems. In particular we study
the Marangoni instability of an infinite
 fluid layer and we show that the numerical results with a 
Chebyshev collocation method are 
highly correspondent to the exact ones. These ideas have been
applied to linear stability analysis of  the B\'enard-Marangoni (BM)
problem in cylindrical geometry and the results obtained have been 
very accurate.

\end{abstract}

\section{Introduction}

The problem of thermoconvective instabilities in fluid layers heated from
below has become a classical subject in fluid mechanics \cite{benard,pearson}. 
It is well known
that two different effects are responsible for the onset of motion when the
temperature difference becomes larger than a certain threshold: gravity and 
capillary forces. When both effects are taken into account the problem is called 
B\'enard-Marangoni (BM) convection \cite{nield}. 
At the beginning theoretical studies 
considered layers of infinite horizontal 
extent without modelling lateral side-wall effects always 
present in experiments. More 
recently experiments have been concerned with confined containers with different 
geometries \cite{kos,ezersky,ordan,marialuisa}, and consecuently theoretical 
and numerical studies have followed the same 
way  \cite{rosenblat1,rosenblat2,dijkstra,dauby1,vrentas,zaman,dauby2}.\par 
Solving numerically hydrodynamical problems of incompressible fluids raises
the question of handling first order derivatives for  pressure in a closed
 container and determining its boundary conditions. In the velocity-pressure 
 formulation \cite{gresho,poz,orszag} first order derivatives
are  avoided by  deriving a Poisson equation for pressure
\cite{gresho}, although the problem of finding
its boundary conditions still remains.  Dirichlet and
Neumann conditions are  studied in detail in Refs. \cite{gresho,orszag} and  
good results
are obtained with Neumann conditions in Ref. \cite{marcus}. In this work, 
the sensitivity of the results, with respect to artificial pressure 
boundary conditions, is reported.
The results we mention refer to rigid boundary conditions 
in the velocity field and purely hydrodynamical problems, but 
the issue of thermoconvection  with free or Marangoni conditions for the velocity 
at the upper surface which we treat here has not been addressed. 
To avoid spurious results 
in the primitive variables formulation the finite difference and finite element methods 
use staggered grids, but this is not necessary with the collocation 
method we propose in this article.

In thermoconvection problems pressure is usually avoided. In Refs. \cite{mer1,mer2}
 the method of potentials of velocity   is used
to eliminate the variable of pressure from the equations. This 
technique raises the order of
the differential equations and additional boundary conditions may be
 required. This is particularly troublesome in cylindrical coordinates
 where high order derivatives cause awkward difficulties.
The spectral method used in Ref. \cite{dauby2} allows 
the removal of pressure in the primitive variables formulation, however 
the collocation method that we study here  is easier to implement. 

In this paper we present several pressure boundary conditions that allow 
us to solve thermoconvective problems in the primitive variables formulation. This is shown 
to be very useful when dealing with cylindrical coordinates.
In the second section we perform the linear stability
analysis of the Marangoni problem in an infinite horizontal fluid layer. 
Since it has an exact solution it is used as a test problem 
to be compared
with the numerical results obtained in the primitive variables formulation and 
several pressure boundary conditions. We use 
a Chebyshev collocation method and sensivity to the pressure 
boundary conditions is observed only
in the  convergence rate of the method. 
In the third section we study the BM problem in cylindrical geometry 
in the primitive variables formulation with the boundary conditions that 
provide better convergence for the test problem. We explain 
the methodology in detail. We study the convergence of the numerical results and
compare them  with the previous bibliography on the subject. 
In the fourth section conclusions are presented.

\section{Test problem}

The physical situation we consider  corresponds to an infinite layer of a fluid
with a free top surface. The domain is:
\begin{eqnarray}
\Omega = \left\{(x,y,z) \in \R^3/ (x,y) \in \R^2,\;0 < z < 1 \right\}.
\end{eqnarray}
In the vertical direction a temperature gradient is imposed.
It is well known that motion sets in after the vertical
temperature gradient has reached a critical value.
The only mechanism that we consider in developing convection is the Marangoni effect,
which consist of the variation of the surface tension with temperature.
In the reference state there is no motion and heat propagates by conduction only.
To study the linear stability of this state we must solve the linearized
equations for the perturbations from the reference state.
As there are no boundary conditions in
the horizontal directions and the basic state is invariant
under translations in those directions it is possible
to take the Fourier modes in those coordinates.
\begin{eqnarray}
{\bf u}' = {\bf u}(z) e^{i {\bf k} {\bf x}} + cc. ,\;\;
\theta' = \theta(z) e^{i {\bf k} {\bf x}} + cc. ,\;\;
p' = p(z) e^{i {\bf k} {\bf x}} + cc.
\end{eqnarray}
The linear problem to be solved is
\begin{eqnarray}
\nabla \cdot{\bf u} &=& 0, \label{eq:continuidad} \\
Pr^{-1}\partial_t {\bf u} &=& - \nabla p +
\Delta {\bf u}, \label{eq:ns}\\
\partial_t \theta &=& u_z + \Delta \theta, \label{eq:calor}
\end{eqnarray}
where ${\bf u}=(u_x,u_y,u_z)$ is the velocity field of the fluid, $\theta$
is the temperature,
$p$ the pressure, $\nabla= (i k_x, i k_y, \partial_z) $, 
$\Delta =
( \partial^2_{z}- k^2)$, $ k = |{\bf k}|$,
$Pr$ is the Prandtl number which is considered infinite.
The boundary conditions are:
\begin{eqnarray}
{\bf u}|_{z=0} &=& 0,\;\theta|_{z=0} = 0, \label{eq:cc1} \\
(\partial_z u_x &+& M i k_x \theta)|_{z=1} = 0,\;\;\;
(\partial_z u_y + M i k_y \theta)|_{z=1} = 0, \label{eq:cc2} \\
 u_z|_{z=1} &=& 0,\;\;\;(\partial_z \theta + B\theta)|_{z=1}=0,
 \label{eq:cc3}
\end{eqnarray}
where $B$ is the Biot number and $M$ is the
Marangoni number which takes into account the variation of the surface 
tension with temperature. 
\subsection{Formulation without pressure}
Classically, pressure is eliminated and the following
set of equations and boundary conditions are obtained,
\begin{eqnarray}
(D^2-k^2)^2 u_z  &=& 0, \label{eq:31}\\
 u_z + (D^2-k^2) \theta & = & 0, \label{eq:32}\\
u_z|_{z=0} = D u_z|_{z=0} = \theta|_{z=0} &=& 0, \label{eq:33}\\
u_z|_{z=1} = (D^2 u_z + k^2 M \theta)|_{z=1} = (D \theta + B
 \theta)|_{z=1} &=&
0, \label{eq:34}
\end{eqnarray}
where $D = \partial/\partial z$. 
This formulation will be called P1 in the following text. 
It can be easily shown that the $u_z$ and $\theta$ solutions are
\begin{eqnarray}
u_z (z) & = & a_1 \,e^{k\,z} + a_2 \,z\,e^{k\,z} + 
a_3 \,e^{ - k\,z} + a_4 \,z\,e^{ - k\,z} \\
\theta (z) &=& b_1 \,e^{k\,z} + b_2 \,e^{- k
\,z} - {\displaystyle \frac {1}{2}} \,{\displaystyle \frac {
a_1 \,z\,e^{k\,z}}{k}}  + {\displaystyle \frac {1}{2}} \,
{\displaystyle \frac {a_3 \,z\,e^{- k\,z}}{k}}  +  \nonumber \\
&& a_2
\,\left( {\displaystyle \frac {1}{4}} \,{\displaystyle \frac {z\,e^{k
\,z}}{k^{2}}}  - {\displaystyle \frac {1}{4}} \,{\displaystyle 
\frac {z^{2}\,e^{k\,z}}{k}} \right) \mbox{} + a_4 \,\left( {\displaystyle \frac {1}{4}} \,
{\displaystyle \frac {z\,e^{ - k\,z}}{k^{2}}}  + 
{\displaystyle \frac {1}{4}} \,{\displaystyle \frac {z^{2}\,e^{
 - k\,z}}{k}} \right)
%\mbox{\hspace{259pt}}
\end{eqnarray}
Introducing them into the equations (\ref{eq:31})-(\ref{eq:34})
we get a linear system of equations for the unknowns $a_i,\,\,b_i$.
The solvability condition for the system gives the following 
dependence of the critical Marangoni number on $k$,
\begin{eqnarray*}
M &=& \frac{8 k B (\alpha + 4 \alpha k + \alpha^2 - 4 \alpha^2 k - 1 -  \alpha^3) + 8k^2( - 1 - \alpha  + 4 k \alpha + \alpha^2  + 4 k \alpha^2  + \alpha^3)}{-3 \alpha^2+ 4 k^3 \alpha + 3 \alpha + 4 k^3 \alpha^2-1+\alpha^3} 
%\lefteqn{{\it M} = 8k({ B}\,\alpha(k) + 4\,{ B}\,
%\alpha(k)\,k + \alpha(k)^{2}\,{ B} - k - \alpha(k)\,k + 4\,k
%^{2}\,\alpha(k) + \alpha(k)^{2}\,k + 4\,k^{2}\,\alpha(k)^{2}} \\
% & & \mbox{} - 4\,{ B}\,\alpha(k)^{2}\,k - { B} - { B
%}\,\alpha(k)^{3} + k\,\alpha(k)^{3})/ \\
% & &  (- 3\,\alpha(k)^{2} + 4\,k^{3}\,\alpha(k) + 3\,\alpha(k)+ 4\,k^{3}\,\alpha(k)^{2} - 1 + 
%\alpha(k)^{3}) \mbox{\hspace{270pt}},
\end{eqnarray*}
where $\alpha$ is a function of $k$, $\alpha(k) = {\rm cosh}(2\,k) - {\rm sinh}(2\,k)$. 

\subsection{Formulation in primitive variables}
We consider the primitive variables formulation
with two different boundary conditions for pressure.
\begin{enumerate}
\item Continuity equation at the boundaries
\begin{eqnarray}
( \nabla \cdot {\bf u})|_{z=0,\,1} = 0 , \label{eq:cc4}
\end{eqnarray}
these boundary conditions together with equations (\ref{eq:continuidad})-(\ref{eq:calor}) and 
boundary conditions (\ref{eq:cc1})-(\ref{eq:cc3}) define the problem P2.
\item Normal component of the Navier-Stokes equations at the top boundary 
and continuity equation at the bottom,
\begin{eqnarray}
(- \partial_z p + \Delta u_z)|_{z=1} &=& 0, \label{eq:cc51}\\
( \nabla \cdot {\bf u})|_{z=0} &=& 0 , \label{eq:cc52}
\end{eqnarray}
these boundary conditions together with equations (\ref{eq:continuidad})-(\ref{eq:calor}) and 
boundary conditions (\ref{eq:cc1})-(\ref{eq:cc3}) define the problem P3.
\end{enumerate}

\subsection{Numerical method}

The exact formulation (P1) together with formulations P2 and P3
are solved numerically with a Chebyshev collocation method. 
The eigenfunctions are approximated by Chebyshev polynomial expansions 
in the $z$ direction.
After changing the $z$ coordinate
 to transform the $[0,1]$ interval into $[-1,1]$ we use the expansions,
\begin{eqnarray}
\displaystyle {\bf u} = \sum_{n=0}^{N-1} {\bf a}_{n} {\rm T}_n (z),
\;\;\theta = \sum_{n=0}^{N-1} b_{n} {\rm T}_n (z),\;\;
p = \sum_{n=0}^{N-1} c_n {\rm T}_n (z),
\end{eqnarray}
which are introduced into the equations.\par
For problem P1 the system and boundary conditions are
evaluated at the collocation points, 
 \begin{eqnarray}
 z_i&=&\cos \left( \left(\frac{i-1}{N-1}-1 \right) \pi \right)
\end{eqnarray}   
In particular, (\ref{eq:31}) is evaluated
at  nodes $i=3,...,N-2$; (\ref{eq:32}) at $i=2,...,N-1$;
the first boundary conditions (\ref{eq:33}) at $i=N$ and the
second ones (\ref{eq:34}) at $i=1$. 
We obtain $2 \times N$ unknowns and $2 \times N$ equations.
For problems P2 and P3:  (\ref{eq:continuidad})
at the nodes $i=2,...,N-1$; (\ref{eq:ns}) at $i=2,...,N-1$; (\ref{eq:calor}) at
$i=2,...,N-1$; the first boundary conditions (\ref{eq:cc1}) at $i=N$;
(\ref{eq:cc2}) and (\ref{eq:cc3}) at $i=1$.
We obtain $5 \times N$ unknowns and $5 \times (N-2)+8$ equations.
To complete the system we evaluate Eq. (\ref{eq:cc4}) at
$i=1,\;N$ for problem P2, and Eq. (\ref{eq:cc51}) at 
$i=1$ and (\ref{eq:cc52}) at
$i=N$ for problem P3. If the coefficients of the unknowns which
form the matrices of the resulting linear system $A$ and $B$
satisfy det$(A- \lambda B)=0$,
a nontrivial solution of the linear homogeneous
system exists. This condition generates a dispersion
relation $\lambda \equiv \lambda(k,M,B)$, 
equivalent to a direct calculation of the eigenvalues from the
system $AX=\lambda BX$,
where $X$ is the vector which contains the
unknowns. When $\lambda$ becomes positive the basic state is
unstable. 
 In the critical situation  $\lambda=0$ with no
imaginary part. In this case the dispersion relation can be written as
$M \equiv M(k,B)$.  We have calculated different critical Marangoni numbers
for several $k$ and $B$ values.

\subsection{Convergence of the numerical method}
To carry out a test of the convergence of the method we
compare the differences in the thresholds of $M_c$ for different orders
of expansions. In table I the thresholds for $k=10$ and $B=10$
are shown for seven odd consecutive 
expansions. The thresholds converge
to the exact value of $1600.01$ and the difference between consecutive 
expansions  tends to zero as $N$ increases. Results on convergence
improves greatly when $k$ and $B$ are smaller, i.e. for the
critical value of $k$, $k_c$ as table II shows, for this reason  there and in table III
we use $7$ polynomials expansions. \par
\subsection{Discussion}

As the exact solution is known, it is 
straightforward to compare the numerical methods. 
In table II 
we show the critical Marangoni and wave numbers in each formulation for different 
Biot numbers. 
We calculate differences between each formulation and the 
exact solution. The mean relative error for P1 and P2 is 
$10^{-2}$, whereas the same for P3 is 
$2 \cdot 10^{-3}$. This indicates P3 as the best approach. 
In table III we show the critical Marangoni number for 
$k=10$ in each formulation and different Biot numbers. 
Although convergence has not 
been reached (see table I) the result is significant  because 
while the mean relative error for P1 is $0.5$ and for 
P2 it is $0.4$,  for P3 it is $0.005$. We confirm again 
that  the third formulation converges better. 
Differences in accuracy
between P1 and P3 could be due to the fact that P1 is a higher order
problem. On the other hand P2 diminishes its efficiency
with respect to P3 because of the boundary conditions considered. Although
 sensitivity
of the numerical method to pressure boundary conditions has been addressed \cite{marcus,mer1}
it  seems here that it is not much more important 
than other effects such as  increasing the  order of 
the derivatives in the equations. Moreover we have shown 
that considering
appropriate boundary conditions for pressure leads
to very good numerical performance.\par
\section{B\'enard-Marangoni in cylindrical geometry}
\subsection{Formulation of the problem}
The physical situation corresponds to a fluid layer of thickness 
$d$ filling a cylinder with radius $l$. 
The surface tension at the upper free surface is temperature 
dependent and the fluid is heated from below. Motion sets in 
when the vertical temperature gradient has exceeded a critical 
value.\par
In the basic state, there is no motion and heat propagates 
by conduction. In the linear Boussinesq 
equations, the field can be expanded in the
azimuthal variable $\phi$ in Fourier modes as follows,
\begin{eqnarray}
{\bf u}(t,r,\phi,z) &=& {\bf u}_m(t,r,z) e^{i m \phi},\\
\theta(t,r,\phi,z) &=& \theta_m (t,r,z) e^{i m \phi},\\
p(t,r,\phi,z) &=& p_m(t,r,z) e^{i m \phi}, \label{eq:um}
\end{eqnarray}
where $\theta$ and ${\bf u}= u {\bf e}_r + v {\bf e}_{\phi} + w {\bf e}_z$ are 
the infinitesimal temperature and velocity perturbations with respect 
to the conductive solution, $p$ is the pressure perturbation, and $(r,\phi,z)$ 
are the polar coordinates. 
If  space, time, velocity, temperature and pressure fields are
respectively divided by the constants $d$, $d^2/\kappa$, $\kappa/d^2$, 
$\Delta T$ and $ \rho_0 \kappa \nu/d^2$, the equations in
dimensionless form are obtained. Here $\kappa$ is the thermal diffusivity, 
$\rho_0$ is the mean density, $\eta$ is the dynamic viscosity (related
to the kinematic viscosity through the expression, $\nu=\eta /\rho_0$)
and $\Delta T$ is the (conductive) temperature drop between the bottom and the top layers. The linearized general equations for the 
perturbations of index $m$ are easily shown to be
\begin{eqnarray}
Pr^{-1} \frac{\partial u}{\partial t} &=& - \frac{\partial p}{\partial r} + 
\nabla_{m^2+1}^2 u  - \frac{2 i m}{r^2} v, \label{eq:1} \\
Pr^{-1} \frac{\partial v}{\partial t} &=& - \frac{i m}{r} p + 
\nabla_{m^2+1}^2 v  + \frac{2 i m}{r^2} u, \label{eq:2} \\
Pr^{-1} \frac{\partial w}{\partial t} &=& -  \frac{\partial p}{\partial z} + 
\nabla_{m^2}^2 w + R \theta, \label{eq:3} \\
\frac{\partial \theta}{\partial t} &=& w + \nabla_{m^2}^2 \theta, \label{eq:4} \\
0 &=& \frac{1}{r} \frac{\partial (r u)}{\partial r} + \frac{i m}{r} v + 
\frac{\partial w}{\partial z} , \label{eq:5} 
\end{eqnarray}
where the subindex $m$ has been cancelled for simplicity. 
The Rayleigh number is defined by $R= \alpha g \Delta T d^3/ \kappa \mu$, 
where $\alpha$ is the coefficient of volume expansion and $g$ is the 
gravity. The Prandtl number is given by $Pr=\nu/ \kappa$. $\nabla_n^2$ is 
defined by $\nabla^2_n =  r^{-1} \partial / \partial r ( r \partial/ \partial r) -
n r^{-2} + \partial^2 / \partial z^2$. \par
The boundary conditions are the following. The bottom of the box is 
rigid and assumed to be perfectly heat conducting so that
\begin{eqnarray}
{\bf u} =0, \; \theta = 0,\;\frac{1}{r} \frac{\partial (r u)}{\partial r} + \frac{i m}{r} v + 
\frac{\partial w}{\partial z}=0 \;\; \mbox{on} \;\;z=0. \label{eq:bcz0}
\end{eqnarray}
The upper surface of the fluid is assumed to  be plane, nondeformable and free where
surface tension effects are taken into account. 
We also assume that at the top, heat is transferred from the liquid to the 
ambient gas according to Newton's law of cooling, which results in a 
Biot condition for the temperature perturbations.  For  pressure the 
normal projection 
of the Navier-Stokes equations into the top plane is considered 
as a boundary condition. 
The mathematical expressions of the 
boundary conditions at the upper surface are then,
\begin{eqnarray}
&& w =0,\; \frac{\partial \theta}{\partial z} + B \theta =0,\; 
\frac{\partial u}{\partial z} + M \frac{\partial \theta}{\partial r} = 0,\;
\frac{\partial v}{\partial z} + M \frac{i m}{r} \theta = 0, \; \mbox{on}\;\;z=1, \label{eq:bcz11}\\
&& -  \frac{\partial p}{\partial z} + \nabla_{m^2}^2 \, w + R \,\theta =0,
\; \mbox{on}\;\;z=1, \label{eq:bcz12}
\end{eqnarray}
where $B$ is the Biot number and $M = \gamma \Delta T d/\rho_0 \kappa \nu$ is 
the Marangoni number with $\gamma$ the constant rate of
change of surface tension with temperature.\par
The lateral side wall is rigid and we will consider it adiabatically 
insulated. For pressure we consider the normal projection 
of the Navier-Stokes equations into this wall as boundary condition. 
This is written as follows
\begin{eqnarray}
&&{\bf u}= {\bf 0},\;\;\partial \theta /\partial r = 0,
 \;\mbox{on}\;\;r= a, \label{eq:bcr11} \\ 
 &&- \frac{\partial p}{\partial r} + 
\nabla_{m^2+1}^2 u  - \frac{2 i m}{r^2} v =0, \;\mbox{on}\;\;r= a,
 \label{eq:bcr12}
\end{eqnarray}
where $a=l/d$ is the aspect ratio. \par
The use of cylindrical coordinates, which are 
singular at $r=0$, imposes regularity conditions on the 
velocity and temperature fields \cite{canuto}. These 
conditions, which express that 
the unknown fields are single valued at $r=0$, are mathematically 
summarized as follows, 
\begin{eqnarray}
\frac{\partial {\bf x}}{\partial \theta}=0,
\end{eqnarray}
where ${\bf x}$ is any scalar or vectorial field.
For the temperature and velocity perturbations, these conditions 
are written as
\begin{eqnarray}
&& u = \frac{\partial w}{\partial r} = \frac{\partial \theta}{\partial r} = 
\frac{\partial p}{\partial r} = 0, \;\;
m =0, \label{eq:bcr0m0} \\
&& u+ i v = w = \theta = p = 0, \;\; m = 1, \label{eq:bcr0m1} \\
&& u = v = w = \theta = p = 0, \;\; m \neq 0,\, 1 . \label{eq:bcr0}
\end{eqnarray}
\subsection{Numerical method}
The numerical method used is a collocation method which is 
similar to the method used in Ref. \cite{icosa} for rectangular 
containers. Firstly we change the variables $z$ and $r$ to transform the
intervals $[0,1] \times [0,a]$ into $[-1,1] \times [-1,1]$ to use 
the Chebyshev polynomials. 
The fields in equations (\ref{eq:1})-(\ref{eq:5}), denoted generically by {\bf x}, 
are aproximated by Chebyshev expansions,
 \begin{eqnarray}
\displaystyle {\bf x} &=& \sum_{l=0}^{L-1} \sum_{n=0}^{N-1} a^{\bf x}_{ln} {\rm T}_l (r) 
{\rm T}_n (z),
\end{eqnarray} 
which are introduced into the equations. The resulting expressions are   
evaluated at the collocation points $(r_j,z_i)$ defined as follows,
 \begin{eqnarray}
 r_j&=&\cos \left( \left(\frac{j-1}{L-1}-1 \right) \pi \right), \\
 z_i&=&\cos \left( \left(\frac{i-1}{N-1}-1 \right) \pi \right).
\end{eqnarray}   
In the case $m > 1$ the system and boundary conditions are evaluated at the 
following collocation points: Eqs. (\ref{eq:1})-(\ref{eq:5}) at the nodes 
$i=2,...,N-1,\; j=2,...,L-1$; the boundary conditions at $z=-1$ (\ref{eq:bcz0}) at 
$i = 1,\;j=2,...,L-1$; the boundary conditions at $z=1$ (\ref{eq:bcz11}) 
at $i=N,\;j=2,...,L-1$ and (\ref{eq:bcz12}) at $i=N,\;j=2,...,L$; 
the boundaries at $r=-1$ (\ref{eq:bcr0})  at 
$i=1,...,N,\;j=1$; finally  the boundaries at $r=1$ (\ref{eq:bcr11}) at 
$i=1,...,N,\;j=L$ and (\ref{eq:bcr12}) at $i=1,...,N-1,\;j=L$. 
We obtain $5 \times N \times L$ equations and $5 \times N \times L $ unknowns. 
For $m=1$ the evaluation has been done at the same nodes, but as there 
are only four boundary conditions at $r=-1$ we have diminished the 
expansion range of the $v$ field by one order, 
\begin{eqnarray}
\displaystyle v = \sum_{l=0}^{L-2} \sum_{n=0}^{N-1} a^v_{ln} 
{\rm T}_l (r) {\rm T}_n (z).
\end{eqnarray} 
 Therefore we get $4 \times N \times L + N \times (L-1)$
equations with the same amount of unknowns. 
In the case $m=0$ the angular component of the velocity is nulle 
and Eq. (\ref{eq:2}), together with the corresponding boundary 
condition for $v$ in Eqs. (\ref{eq:bcz0})-(\ref{eq:bcr0}) disappear and therefore 
we obtain $4 \times N \times L$ equations with the same amount of 
unknowns. In this case the matrix associated to the 
linear algebraic system is singular, due to the fact that 
pressure is defined up to an additive constant. To fix this 
constant in the node $i=N-2, j=L$ the boundary condition (\ref{eq:bcr12}) 
is replaced by a Dirichlet condition for pressure ($i.e$, $p=0$ at $i=N-2, j=L$). 
The resulting linear systems have been solved with a standard numerical package. 
The eigenfunctions and thresholds of the generalized problem are numerically calculated. 
If the coefficients of the unknowns which 
form the matrices $A$ and $C$ satisfy 
det$(A- \lambda C)=0$, 
a nontrivial solution of the linear homogeneous 
system exists. This condition generates a dispersion 
relation $\lambda \equiv \lambda(R,M,B)$, 
equivalent to a direct calculation of the eigenvalues from the 
system $AX=\lambda CX$, 
where $X$ is the vector which contains the 
unknowns. When $\lambda$ becomes positive the basic state is 
unstable. In the critical situation  $Re(\lambda)=0$ without
imaginary part. In this case the dispersion relation can be written as
$M \equiv M(R,B)$ or $R \equiv R(M,B)$ and critical Marangoni number or
critical Rayleigh number are obtained directly from the eigenvalue problem.

\subsection{Convergence of the numerical method}
To carry out a test of the convergence of the method we
compare the differences in the critical Marangoni number $M_c$ 
to different orders of expansions for different values of 
the aspect ratio $a$. In table IV these thresholds 
are shown for four consecutive
expansions varying the number of polynomials taken in the 
$z$ direction ($N$) and in the $r$ direction ($L$). 
These results allow us to conclude that $N \times L = 9 \times 13$ 
gives very good results for aspect ratios lower than $5$. So, except 
where otherwise stated, all results given below correspond to 
the values $N \times L = 9 \times 13$.
\subsection{Comparison with other theoretical works}
Our main goal is to show how a Chebyshev collocation method 
applied to the primitive variables formulation of 
thermoconvective problems with convenient boundary conditions 
for pressure produces excellent results. For this reason 
we compare our results with the preceding ones on 
the subject.\par
% There are several papers on the same problem, 
% we compare here with those of Vrentas {\it et al.} for 
% $m=0$, Zaman and Narayanan for $m \neq 0$ and Dauby {\it et al.} for 
% any $m$. \par
Dauby {\it et al.} have already considered the comparison between 
their results and the preceding ones. For this reason we are 
focusing on a comparison only with their paper. 
The correspondence between our work and that of Dauby {\it et al.} can be tested 
by comparing our Tables V and VI to their Tables II and III. 
>From this comparison we see that the deviations are less than $0.06 \%$. 
We have checked that the succession of unstable eigenmodes 
when $a$ is increasing  in our  approach  coincides with
those of Dauby {\it et al.}. In their work they discuss 
 differences obtained  with Zaman and Narayanan \cite{zaman}. 
As is shown in figure 1 when $a$ increases all the critical 
values converge to the  same thresholds and it is
 harder  to distinguish among them with a numerical method.
\section{Conclusions}
We have shown that a Chebyshev collocation method 
applied to the primitive variables formulation of 
thermoconvective problems with convenient boundary conditions 
for pressure produces excellent results. 
First, we have solved the Marangoni problem in infinite geometry 
with a Chebyshev collocation method, 
keeping the primitive variables formulation and testing two sets of boundary conditions
for pressure. We have compared these results with the exact ones 
obtained classically by eliminating pressure. 
As there is an exact solution to 
this problem it was straightforward to compare the performance 
of the  numerical results in the different formulations. 
We prove that the results obtained do not depend on those 
boundary conditions. Sensitivity to those boundary conditions 
is only observed in the convergence of the method. 
The best results are obtained 
 in the primitive variables problem, taking as boundary condition 
for pressure the projection of Navier-Stokes equations at the 
top boundary and continuity at the bottom one. 
Second, considering the boundary 
condition with the best convergence we have solved the BM
problem in cylindrical geometry. We have compared with 
previous numerical 
results and the correspondence is very high. 
We conclude that the primitive equations 
with the appropriate boundary conditions for pressure 
of these thermoconvective problems 
are accurately solved with our collocation method.
\subsection*{Acknowledgments}
We are deeply grateful to Professor  
Guo Ben-Yu and to D. Valladares, N. L\'opez, M. Net and 
F. Marqu\'es for our rewarding discussions.  
This work was partially supported by the  
DGICYT (Spanish Government) under Grant No. PB96-0534 
and by the University of Castilla-La~Mancha.

\newpage
\noindent
\centerline{\bf Table I}
%\begin{table}
\begin{center}
\begin{tabular}{cccc} \hline \hline
$N$ & P1 & P2 & P3 \\
\hline 
$17$  & $1600.01$ & $1600.01$ & $1600.01$  \\
$15$  & $1599.99$ & $1599.91$ & $1600.00$ \\
$13$  & $1599.42$ & $1597.95$ & $1599.66$  \\
$11$  & $1590.73$ & $1572.61$ & $1595.49$  \\
$9$   & $1500.53$ & $1382.53$ & $1572.57$  \\
$7$   & $972.58$  & $776.62$  & $1622.50$  \\
$5$   & $231.32$  & $201.84$  & $2817.82$  \\
\hline \hline
\end{tabular}
\end{center}
%\caption{}
%\end{table}
\bigskip
\noindent
%\begin{table}
\centerline{\bf Table II}
\begin{center}
\begin{tabular}{ccccccccc} \hline \hline
$B$ & \multicolumn{2}{c}{P1 (numerical)} & 
\multicolumn{2}{c}{P1 (exact)}
 & 
\multicolumn{2}{c}{P2}  
 &\multicolumn{2}{c}{P3}  
\\ 
\hline \
      & $M_c$   & $k_c$  & $M_c$    & $k_c$ & $M_c$   & $k_c$ & $M_c$   & $k_c$   \\
\hline \\
$0.1$ & 83.20 & 2.04 & 83.43 & 2.03  & 82.93 & 2.05  & 83.31 & 2.03  \\
$1$ & 115.69 & 2.26 & 116.13 & 2.25  & 115.14 & 2.28 & 115.92 & 2.25  \\
$10$ & 410.45 & 2.77 & 413.44 & 2.74 & 406.27 & 2.82 & 412.20 & 2.75  \\
\hline \hline
\end{tabular}
\end{center}
%\caption{}
%\end{table}
\bigskip
\noindent
%\begin{table}
\centerline{\bf Table III}
\begin{center}
\begin{tabular}{ccccc} \hline \hline
$B$ & P1 (numerical) & P1 (exact) & P2 & P3 \\
\hline 
$0.1$  & $484.39$ & $808.01$  & $386.80$ & $808.08$  \\
$1$    & $528.77$ & $880.01$  & $422.23$ & $882.12$  \\
$10$   & $972.58$ & $1600.01$ & $776.62$ & $1622.50$  \\
\hline \hline
\end{tabular}
\end{center}
%\caption{}
%\end{table}
\centerline{\bf Table IV}
\begin{center}
\begin{tabular}{ccccc} \hline \hline
$m$ & $5 \times 9$ & $7 \times 11$ & $9 \times 13$ &  $11 \times 15$ \\
\hline 
$m=0$  & $150.705$ & $153.689$  & $154.064$ &   $154.062$ \\
$m=1$   & $148.156$ & $152.558$  & $152.945$ &  $152.949$ \\
$m=2$   & $150.946$ & $153.734$ & $154.154$ &  $154.154$ \\
$m=3$  &  $148.526$ & $152.899$ & $153.255$ &  $153.256$ \\
$m=4$ & $149.400$  & $153.542$ & $153.979$ &  $153.987$ \\
\hline \hline
\end{tabular}
\end{center}
\centerline{\bf Table V}
\begin{center}
\begin{tabular}{ccccccc} \hline \hline
& \multicolumn{2}{c}{$a=1$ }& 
\multicolumn{2}{c}{$a=2$}& 
\multicolumn{2}{c}{$a=4$} \\  
\hline
 & $M_c$ & $M_c$ &$M_c$ &$M_c$ &$M_c$ &$M_c$ \\
$B(R=0)$ & (H. and M.) & (D. {\it et al.}) & (H. and M.) & (D. {\it et al.}) 
& (H. and M.) & (D. {\it et al.}) \\
\hline
$0.01$  & $164.65$ & $164.55$	&	$84.640$ & 	$84.638$	&	$82.485$ & $82.486$ \\
$0.1$   & $168.44$ & $168.34$	&	$88.346$ & 	$88.344$	&	$86.263$ & $86.263$ \\
$1$     & $206.29$ & $206.16$	&	$125.011$ & $125.000$	&	$120.624$& $120.630$ \\
 & $R_c$ & $R_c$ & $R_c$ & $R_c$ & $R_c$ & $R_c$ \\
$B(M=0)$ & (H. and M.) & (D. {\it et al.}) & (H. and M.) & (D. {\it et al.}) 
& (H. and M.) & (D. {\it et al.}) \\
\hline
$0.01$  & $1419.30$ & $1419.47$ & $712.542$ & $712.667$ & $695.543$ & $695.668$ \\
$0.1$   & $1426.06$ & $1426.24$ & $726.574$ & $726.704$ & $709.326$ & $709.457$ \\
$1$     & $1481.89$ & $1482.12$ & $835.897$ & $836.072$ & $799.522$ & $799.735$ \\
\hline \hline
\end{tabular}
\end{center}
\centerline{Table VI}
\begin{center}
\begin{tabular}{ccccc} \hline \hline
 & $a=1$ & $a=2$ & $a=4$ & $a=8$ \\
\hline 
$m=0$  & $163.676$ & $80.878$ & $78.777$ & $76.179$ \\
$m=1$  & $108.383$ & $91.254$ & $77.864$ & $76.448$ \\
$m=2$  & $158.994$ & $98.407$ & $79.699$ & $76.203$ \\
$m=3$  & $255.885$ & $99.955$ & $78.095$ & $76.487$ \\
\hline \hline
\end{tabular}
\end{center}
\newpage
\noindent
{\bf TABLE CAPTIONS}\par
\bigskip
\noindent
{\bf Table I}\par
\noindent
Critical Marangoni number $M_c$ for different orders of 
expansions and $B=10,\;k=10$.\par
\bigskip
\noindent
{\bf Table II}\par
\noindent
Critical Marangoni number $M_c$ and wave number $k_c$ for 
different values of the Biot number $B$.\par
\bigskip
\noindent
{\bf Table III}\par
\noindent
Critical Marangoni number $M_c$ for 
different values of the Biot number $B$ and $k=10$.\par
\bigskip
\noindent
{\bf Table IV}\par
\noindent
Critical Marangoni number $M_c$ for different orders of 
expansions and $B=2,\;a=5$.\par
\bigskip
\noindent
{\bf Table V}\par
\noindent
Comparison between our results (H. and M.) and those of Dauby {\it et al.}
\cite{dauby2}. The first lines give the critical Marangoni number when 
$R=0$ and for different Biot numbers; the bottom of the table 
presents the critical Rayleigh number for $M=0$.\par
\bigskip
\noindent
{\bf Table VI}\par
\noindent
Critical Marangoni numbers $M_c$ for different azimuthal 
wave number $m$. The Rayleigh and Biot numbers are 
100 and 0.2, respectively.\par
\bigskip
\noindent
{\bf FIGURE CAPTIONS}\par
\bigskip
\noindent
{\bf Figure 1}\par
\noindent
Critical Marangoni number $M_c$ as a function of the aspect ratio 
$a$. Curves corresponding to different azimuthal wave numbers 
$m$ are represented. The Rayleigh and Biot numbers are zero.
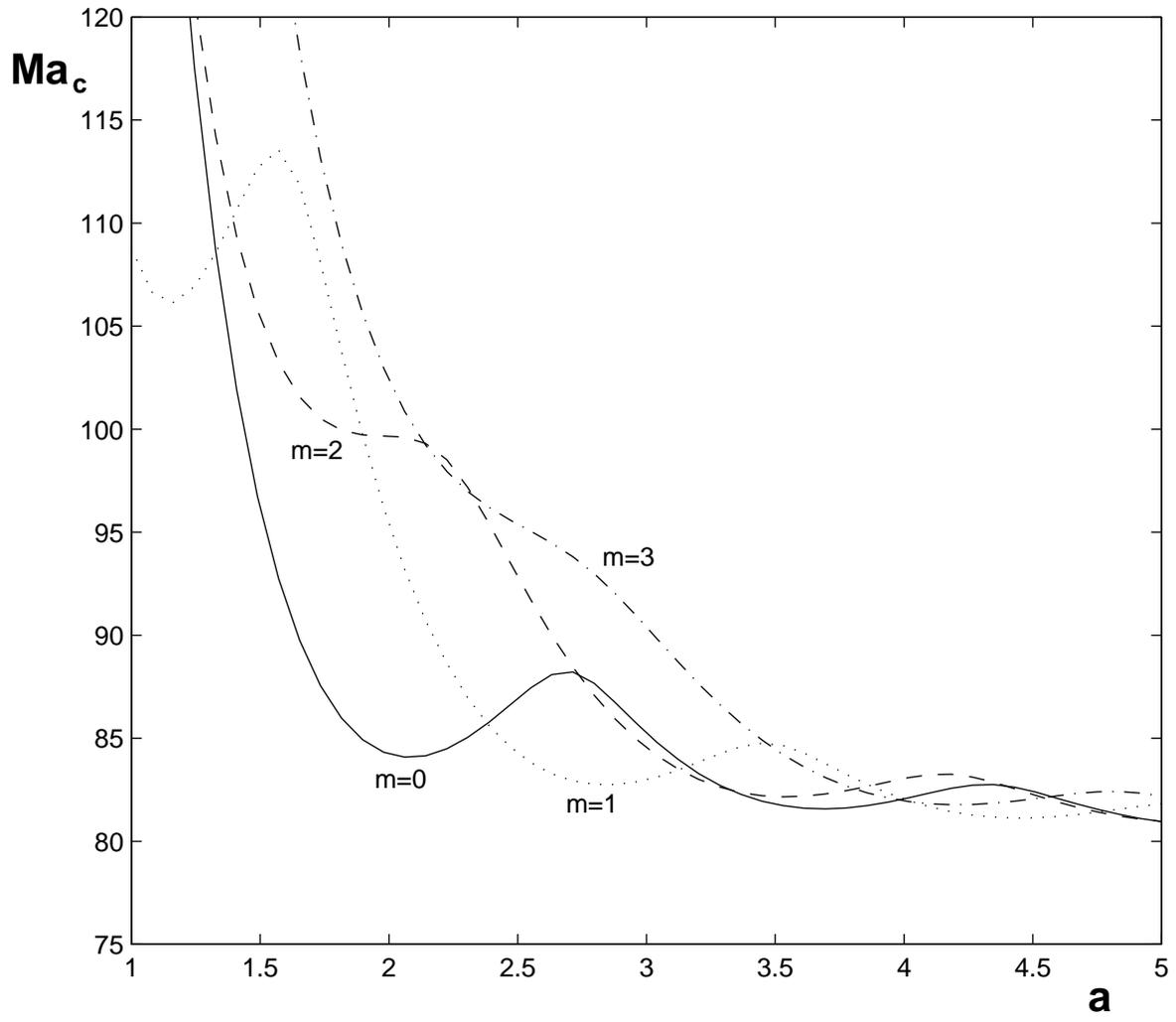
\end{document}